\documentclass[a4paper]{article}

\usepackage{amsmath}
\usepackage{amssymb}
\usepackage{graphicx}
\usepackage{amsthm}
\usepackage{mathtools}
\usepackage{graphicx}
\usepackage{mathrsfs}
\usepackage{amsfonts}
\usepackage{enumerate}
\usepackage{latexsym, euscript, epic, eepic, color}
\usepackage{hyperref}
\usepackage{authblk}
\usepackage{blindtext}

\newtheorem{theorem}{Theorem}[section]

\newtheorem{proposition}{Proposition}[section]

\newtheorem{corollary}{Corollary}[section]

\numberwithin{equation}{section}

\def\brackets#1{\left( #1 \right)}

\def\absolutevalues#1{\left| #1 \right|}

\title{On H\"older maps and prime gaps}

\author{Haipeng Chen$^1$ \  and \ Jonathan M. Fraser$^2$}
\affil{$^1$College of Big Data and Internet, \\ Shenzhen Technology University, Shenzhen, China.\\ hpchen0703@foxmail.com \\
$^2$Mathematical Institute, University of St Andrews, UK \\ jmf32@st-andrews.ac.uk}
\begin{document}

\pagenumbering{arabic}

\maketitle

\begin{abstract}
Let $p_n$ denote the $n$th prime, and consider the function $1/n \mapsto 1/p_n$ which maps the reciprocals of the positive integers bijectively to the reciprocals  of the primes.   We show that H\"older continuity of this function is equivalent to a parametrised family of Cram\'er type estimates on the gaps between successive primes.  Here the parametrisation comes from the H\"older exponent. In particular, we show that Cram\'er's conjecture is equivalent to the map $1/n \mapsto 1/p_n$ being  Lipschitz. On the other hand, we show that the inverse map $1/p_n \mapsto 1/n$ is H\"older of all orders but not Lipschitz and this is independent of Cram\'er's conjecture.
\end{abstract}

\section{Cram\'er's conjecture, H\"older maps,  and our main result}

Understanding the asymptotic  properties of the primes is a fundamental and multifaceted problem in number theory. Let $\{p_n \}_{n = 1}^{\infty}$  denote  the set of primes where $p_n$ is the $n$th prime number and $p_{n+1} > p_n$ for all $n$. Recall the  \emph{Prime Number Theorem} (PNT), which describes the asymptotic growth rate of $p_n$, and  \emph{Rosser's Theorem}, which bounds the $n$th prime by
\begin{equation*}\label{ROSSER.THEOREM.}
n\brackets{\log n + \log \log n -\frac{3}{2}} \leq p_n \leq n\brackets{\log n + \log \log n - \frac{1}{2}}
\end{equation*}
for all $n \geq 20$. For   further discussion of  these results see \cite{HW2008, RS1962} and references therein.

A related problem is to consider the gaps between successive primes, see \cite{FGKMT2018,M2015,M2016,Z2014}.  \emph{Cram\'er's conjecture} asserts that there should exist a constant $C>0$ such that 
\[
p_{n+1} - p_n \leq C (\log p_n )^2
\]
for all $n \geq 1$.   In particular, using Rosser's theorem,  Cram\'er's conjecture gives
\[
p_{n+1} - p_n \leq C' (\log n )^2
\]
for all $n \geq 1$ and a different constant $C'$. The main objective of this paper is to connect Cram\'er's conjecture to a problem concerning H\"older exponents of the natural map between the reciprocals of the positive integers and the reciprocals of the primes.  This approach is motivated by various problems in metric geometry where one tries to understand a given metric space by identifying  those spaces which are in the same  bi-Lipschitz  equivalence class.  For example, bi-Lipschitz equivalence implies coincidence of  familiar notions of fractal dimension such as Hausdorff, box and Assouad dimension. To this end we consider the bi-H\"older continuity of the map $1/n \mapsto 1/p_n$  mapping the reciprocals of the positive integers bijectively to  the reciprocals of the primes.  Our first result proves this map has a H\"older inverse of all orders.  Recall that a map $f: X \to Y$ is \emph{H\"older} (of order $\alpha \in (0,1)$) if there exists a constant $c \geq 1$ such that
\[
|f(x)-f(y)| \leq c|x-y|^\alpha
\]
for all $x,y \in X$.  We assume here that $X$ and $Y$ are bounded subsets of Euclidean space, but this is not necessary in general.  The map $f$ is \emph{bi-H\"older} if it is H\"older and has a H\"older inverse and \emph{Lipschitz} if it satisfies the  H\"older condition but with the optimal order $\alpha=1$. 

\begin{theorem}\label{LOWER.BOUND.}
For all  $\varepsilon > 0$ there exists an  integer $N(\varepsilon)$ such that, for all $m > n > N(\varepsilon)$, we have
\[
\absolutevalues{\frac{1}{p_n} - \frac{1}{p_m}} \geq \frac{1}{6} \cdot  \absolutevalues{\frac{1}{n} - \frac{1}{m}}^{1+\varepsilon}.
\]
\end{theorem}

We note that Theorem \ref{LOWER.BOUND.} is sharp in the sense that it cannot be `upgraded' to a Lipschitz bound. For example, results on bounded gaps between primes, e.g.  \cite{Z2014},  show that
\[
\liminf_{n \to \infty} (p_{n+1}-p_n) < \infty.
\]
In particular, applying this result together with the PNT yields
\[
\liminf_{n \to \infty}\frac{\frac{1}{p_n} - \frac{1}{p_{n+1}}}{\frac{1}{n} - \frac{1}{n+1}} = \liminf_{n \to \infty} \frac{n(n+1)}{p_np_{n+1}} \cdot (p_{n+1}-p_n) =0.
\]
This proves that the map $1/p_n \mapsto 1/n$ is not Lipschitz.

H\"older continuity of the forward map is more subtle.  Our next result shows that if $m$ and $n$ are `sufficiently separated', then a bi-Lipschitz estimate can be derived, up to a logarithmic error.

\begin{theorem}\label{HOLDER.ESTIMATE.LARGE.SCALE.}
For all  $0 < \varepsilon < 1$ there exist an integer $N(\varepsilon)$ such that for all $n \geq N(\varepsilon)$ and $m > \frac{1 +  \varepsilon}{1-\varepsilon}n$, we have 
\[
\frac{1}{(1+\varepsilon)^2} \cdot \brackets{\frac{1}{n} - \frac{1}{m}} \frac{1}{\log m} \leq \frac{1}{p_n} - \frac{1}{p_m} \leq \brackets{1+\varepsilon} \cdot  \brackets{\frac{1}{n} - \frac{1}{m}} \frac{1}{\log n}.
\]
\end{theorem}

It follows immediately  from Theorem \ref{HOLDER.ESTIMATE.LARGE.SCALE.} that the forward map is actually Lipschitz continuous  in the range $m \geq 2n$  for sufficiently large $n$.
\begin{corollary}\label{cor}
For all  $0 < \varepsilon < 1$ there exists an integer  $N(\varepsilon)$ such that for all $n \geq N(\varepsilon)$ and $m > \frac{1 +  \varepsilon}{1-\varepsilon}n$, we have 
\[
\frac{1}{p_n} - \frac{1}{p_m} \leq \frac{1 + \varepsilon}{\log N(\varepsilon)} \brackets{\frac{1}{n} - \frac{1}{m}}.
\]
\end{corollary}

 H\"older continuity of the forward map over the full range is related to a parametrised family of Cram\'er type bounds on  prime gaps.  This is the content of Theorems \ref{UPPER.BOUND.1}-\ref{UPPER.BOUND.2}.

\begin{theorem}\label{UPPER.BOUND.1}
Suppose that for  $\varepsilon \geq 0$, there exists a constant $c(\varepsilon)$ and an  integer $N(\varepsilon) \geq 20$ such that, for all $n > N(\varepsilon)$,
\[
\absolutevalues{\frac{1}{p_n} - \frac{1}{p_{n+1}}} \leq c(\varepsilon) \absolutevalues{\frac{1}{n} - \frac{1}{n+1}}^{1-\varepsilon}.
\]
Then, for all sufficiently large $n$, we have
\[
p_{n+1} - p_n \leq 2c(\varepsilon) (n^\varepsilon \log n )^2.
\]
\end{theorem}

Conversely, the weaker forms of Cram\'er's conjecture imply  the forward map is H\"older of all orders. 

\begin{theorem}\label{UPPER.BOUND.2}
Suppose that for  $\varepsilon \geq 0$, there exists a  constant $c(\varepsilon)$ and an integer $N(\varepsilon) \geq 20$ such that, for all  $n > N(\varepsilon)$, 
\[
p_{n+1} - p_n \leq c(\varepsilon) (n^\varepsilon \log n )^2.
\]
Then, for all sufficiently large $n$, we have
\[
\absolutevalues{\frac{1}{p_n} - \frac{1}{p_{n+1}}} \leq c(\varepsilon) \absolutevalues{\frac{1}{n} - \frac{1}{n+1}}^{1-\varepsilon}.
\]
Moreover, for all  $  n < m < 2n$ with $n$ sufficiently large, we  have 
\[
\absolutevalues{\frac{1}{p_n} - \frac{1}{p_m}} \leq c(\varepsilon) \absolutevalues{\frac{1}{n} - \frac{1}{m}}^{1-2\varepsilon}.
\]
\end{theorem}

Combining Theorems \ref{UPPER.BOUND.1}, \ref{UPPER.BOUND.2} and Corollary \ref{cor}, we note that the $\varepsilon = 0$ case shows that Cram\'er's conjecture is equivalent to the map $1/n \mapsto 1/p_n$ being  Lipschitz.

In light of Theorem  \ref{UPPER.BOUND.2}, it is natural to ask for which $\varepsilon>0$ are these weak Cram\'er bounds known to hold.  Note that Bertrand's postulate may be regarded as the first step in this line of research, verifying the case $\varepsilon = 1$.  To the best of our knowledge the state of the art here is provided by Baker,   Harman and  Pintz \cite{bhp} who proved that the interval $[n, n+n^{0.525}]$ always contains a prime for sufficiently large $n$.  In particular, combining this with Theorem  \ref{UPPER.BOUND.2}, Corollary \ref{cor} and Theorem \ref{LOWER.BOUND.} yields the following corollary.

\begin{corollary}
For all $\varepsilon>0$, there exists a constant $C = C(\varepsilon) >0$  such that, for all   $m > n \geq 1$,
\[
C^{-1} \absolutevalues{\frac{1}{n} - \frac{1}{m}}^{1+\varepsilon} \leq \absolutevalues{\frac{1}{p_n} - \frac{1}{p_m}} \leq C \absolutevalues{\frac{1}{n} - \frac{1}{m}}^{0.475}.
\]
\end{corollary}

\begin{proof}
The result of Baker,   Harman and  Pintz \cite{bhp} together with the PNT implies that for some constant $C'>0$ we have
\[
p_{n+1}-p_n \leq p_n^{0.525} \leq C'(n^{0.525/2} \log n)^2
\]
for sufficiently large $n$.  Applying Theorem  \ref{UPPER.BOUND.2} (with $\varepsilon = 0.525/2$ and so $1-2 \varepsilon = 0.475$) proves the desired upper bound for sufficiently large $n$ and $n<m<2n$.   Corollary \ref{cor} takes care of the case when $m \geq 2n$.  Theorem \ref{LOWER.BOUND.} provides  the lower bound for sufficiently large $n$.  Finally, the result follows by ensuring $C$ is chosen large enough to also deal with the small $n$.
\end{proof}

 Theorem \ref{LOWER.BOUND.} is proved in Section \ref{LOWER.BOUND.proof}.   Theorem \ref{HOLDER.ESTIMATE.LARGE.SCALE.} is proved in Section \ref{HOLDER.ESTIMATE.LARGE.SCALE.proof}. Finally, Theorems   \ref{UPPER.BOUND.1} and \ref{UPPER.BOUND.2} are proved in Section \ref{UPPER.BOUND.proof}.

\section{Proof of  Theorem \ref{LOWER.BOUND.}:  H\"older continuity of inverse map} \label{LOWER.BOUND.proof}

In this section, we prove Theorem \ref{LOWER.BOUND.}, which uses Rosser's theorem and the convex version of Jensen's inequality.

\begin{proof}[Proof of Theorem \ref{LOWER.BOUND.}]
Fix $\varepsilon > 0$.  For sufficiently large $n$, we have
\[
2 (\log n) \log (n+1) \leq (n(n+1))^\frac{\varepsilon}{2}.
\]
Thus, by Rosser's theorem,  for sufficiently large $n$
\[
\absolutevalues{\frac{1}{p_n} - \frac{1}{p_{n+1}}}  \geq \frac{1}{p_np_{n+1}} \geq  \frac{1}{n(n+1)} \cdot \frac{1}{2(\log n)\log (n+1)} \geq  \brackets{\frac{1}{n(n+1)}}^{1+\frac{\varepsilon}{2}} = \absolutevalues{\frac{1}{n} - \frac{1}{n+1}}^{1+\frac{\varepsilon}{2}}. 
\]
We now consider two cases, assuming that $n$ is sufficiently large for the above to hold.

\noindent {\bf Case 1. $n < m \leq 2n$.} 
It follows from the convex version of  Jensen's inequality that
\begin{align*}
\absolutevalues{\frac{1}{p_n} - \frac{1}{p_m}} & = \absolutevalues{\frac{1}{p_n} - \frac{1}{p_{n+1}}} + \dots + \absolutevalues{\frac{1}{p_{m-1}} - \frac{1}{p_m}}  \\
& \geq \absolutevalues{\frac{1}{n} - \frac{1}{n+1}}^{1+\frac{\varepsilon}{2}} + \dots + \absolutevalues{\frac{1}{m-1} - \frac{1}{m}}^{1+\frac{\varepsilon}{2}} \\
& = \brackets{m-n} \cdot \frac{1}{m-n} \cdot \brackets{\absolutevalues{\frac{1}{n} - \frac{1}{n+1}}^{1+\frac{\varepsilon}{2}} + \dots + \absolutevalues{\frac{1}{m-1} - \frac{1}{m}}^{1+\frac{\varepsilon}{2}}} \\
& \geq \brackets{m-n} \cdot \brackets{\frac{1}{m-n}}^{1+\frac{\varepsilon}{2}} \cdot \brackets{\frac{1}{n} - \frac{1}{m}}^{1+\frac{\varepsilon}{2}} \\
& = \brackets{\frac{1}{m-n}}^{\frac{\varepsilon}{2}}  \cdot  \brackets{\frac{1}{n} - \frac{1}{m}}^{1+\frac{\varepsilon}{2}}.
\end{align*}
Since $n < m \leq 2n$, we have
\[
\brackets{\frac{1}{m-n}}^{\frac{\varepsilon}{2}}  \cdot  \brackets{\frac{1}{n} - \frac{1}{m}}^{1+\frac{\varepsilon}{2}}  = \brackets{\frac{1}{m-n}}^{\frac{\varepsilon}{2}}  \cdot  \brackets{\frac{1}{n} - \frac{1}{m}}^{-\frac{\varepsilon}{2}} \cdot  \brackets{\frac{1}{n} - \frac{1}{m}}^{1+\varepsilon} 
 \geq \brackets{\frac{1}{n} - \frac{1}{m}}^{1+\varepsilon}.
\]
\noindent {\bf Case 2 $m > 2n$.} 
It follows from Rosser's Theorem that for sufficiently large $n$
\[
\absolutevalues{\frac{1}{p_n} - \frac{1}{p_m}}   \geq \absolutevalues{\frac{1}{(3/2) n \log n} - \frac{1}{2n \log 2n}} 
\geq \frac{1}{6} \cdot \frac{1}{n \log n} 
\geq  \frac{1}{6} \cdot \brackets{\frac{1}{n}}^{1+\varepsilon} 
\geq  \frac{1}{6} \cdot \absolutevalues{\frac{1}{n} - \frac{1}{m}}^{1+\varepsilon}.
\]
Taking case 1 and 2 together proves the result. 
\end{proof}

\section{Proof of Theorem \ref{HOLDER.ESTIMATE.LARGE.SCALE.}: Lipschitz continuity of forward map} \label{HOLDER.ESTIMATE.LARGE.SCALE.proof}

We first  provide an estimate for gaps between primes which will be used in the subsequent proof.  We remark that the estimate  is false if $m$ and $n$ are not assumed to be sufficiently separated.  For example, results on bounded gaps between primes, e.g.  \cite{Z2014},  show that
\[
\liminf_{n \to \infty} (p_{n+1}-p_n) < \infty.
\]
In particular, $p_{n+1}-p_n$ cannot be bounded below by any function which grows without bound, such as $\log n$.

\begin{proposition}\label{GAP.PRIME.INTEGER.CASE.}
For all $0<\varepsilon <1$, there exists an integer $N(\varepsilon)$ such that for all $n \geq N(\varepsilon)$ and $\frac{1+\varepsilon}{1-\varepsilon} n < m$, we have
\[
(m-n)\log n \leq p_m - p_n \leq (1+\varepsilon) (m-n) \log m.
\]
\end{proposition}

\begin{proof}
Fix  $0< \varepsilon <1$.  By Rosser's Theorem there exists an integer $N(\varepsilon)$ such that for any $n > N(\varepsilon)$, we have
\[
n\brackets{\log n + (1-\varepsilon)\log \log n} \leq p_n \leq n\brackets{\log n + \log \log n }
\]
and 
\[
1+\log n \leq \brackets{1+\frac{\varepsilon}{4}} \log n, \qquad
\log \log n + \frac{1}{\log n} \leq \frac{\varepsilon}{2} \log n.
\]
Then, for any $m \geq \frac{1+\varepsilon}{1-\varepsilon} n$, we obtain 
\[
\brackets{m\log m - n \log n} + \brackets{(1-\varepsilon) m\log \log m -  n \log \log n}  \leq p_m - p_n 
\]
and
\[  p_m - p_n  \leq \brackets{m\log m - n \log n} + \brackets{m\log \log m - (1-\varepsilon) n \log \log n}.
\]
Consider the upper bound.  Let $c_1(\varepsilon) = \frac{3}{2} -\frac{\varepsilon}{2} > 1$ noting that
\[
m\log \log m - (1-\varepsilon) n \log \log n \leq c_1(\varepsilon) \brackets{m \log \log m - n \log \log n},
\]
and
\[
\frac{c_1(\varepsilon)-1}{c_1(\varepsilon)-1+\varepsilon} = \frac{1-\varepsilon}{1+\varepsilon} \geq \frac{n\log \log n}{m\log \log m}.
\]
The lower bound can be handled similarly.  Let  $c_2(\varepsilon) = \frac{1 - \varepsilon}{2}\in (0,1)$ noting that
\[
(1-\varepsilon) m \log \log m - n\log \log n \geq c_2(\varepsilon) (m\log \log m - n\log \log n), 
\]
and
\[
\frac{1-c_2(\varepsilon)-\varepsilon}{1-c_2(\varepsilon)} = \frac{1-\varepsilon}{1+\varepsilon} \geq \frac{n \log\log n}{m \log\log m}.
\]
Then, by mean value theorem applied to  $x \log x$ and $x \log \log x$, we obtain
\[
(1 + \log n )(m-n) \leq m \log m - n \log n \leq (1 + \log m )(m-n)
\]
and
\[
\brackets{\log \log n + \frac{1}{\log n}}(m-n) \leq m\log \log m-n \log \log n \leq \brackets{\log \log m + \frac{1}{\log m}}(m-n).
\]
Thus for the upper bound, we have
\[
p_m -p_n \leq (m-n)(1 + \log m) + c_1(\varepsilon) (m-n) \brackets{\log \log m + \frac{1}{\log m}} \leq  (1+\varepsilon) \cdot (m-n) \log m,
\]
Similarly, for the lower bound, we have
\[
p_m -p_n \geq (m-n)(1 + \log n) + c_2(\varepsilon) (m-n) \brackets{\log \log n + \frac{1}{\log n}} \geq (m - n)\log n,
\]
completing the proof.
\end{proof}

We are now ready to prove Theorem \ref{HOLDER.ESTIMATE.LARGE.SCALE.}.

\begin{proof}[Proof of Theorem \ref{HOLDER.ESTIMATE.LARGE.SCALE.}]
Fix $\varepsilon > 0$.  By appealing to  Rosser's Theorem  we may choose an integer $N'(\varepsilon) \geq N(\varepsilon)$ where $N(\varepsilon)$ is the constant from Proposition \ref{GAP.PRIME.INTEGER.CASE.} such that, for any $n \geq N'(\varepsilon)$,  $n \log n \leq p_n \leq (1+\varepsilon) n \log n$.

Consider   $n \geq N'(\varepsilon)$ and  $m \geq \frac{1+\varepsilon}{1-\varepsilon} n$. For the upper bound, it follows from Proposition \ref{GAP.PRIME.INTEGER.CASE.} that 
\[
\frac{1}{p_n} - \frac{1}{p_m} \ = \ \frac{p_m-p_n}{p_m p_n}  \leq (1+\varepsilon) \cdot \frac{(m-n) \log m}{mn\log m \log n} \leq \brackets{1+\varepsilon} \brackets{\frac{1}{n}-\frac{1}{m}} \frac{1}{\log n}.
\]
The lower bound  is similar. It also follows from Rosser's Theorem and Proposition \ref{GAP.PRIME.INTEGER.CASE.} that
\[
\frac{1}{p_n} - \frac{1}{p_m} = \frac{p_m-p_n}{p_m p_n} 
\geq \frac{1}{(1+\varepsilon)^2} \cdot \frac{(m-n)\log n}{mn\log m \log n} =  \frac{1}{(1+\varepsilon)^2} \cdot \brackets{\frac{1}{n} - \frac{1}{m}}\frac{1}{\log m}
\]
completing the proof.
\end{proof}

\section{Proofs of Theorems \ref{UPPER.BOUND.1}-\ref{UPPER.BOUND.2}:  H\"older continuity of forward map and Cram\'er type estimates} \label{UPPER.BOUND.proof}

Theorem \ref{UPPER.BOUND.1} shows that H\"older continuity of the  forward map $1/n \mapsto 1/p_n$ implies a weak form of Cram\'er's conjecture.  This follows easily from Rosser's theorem.
\begin{proof}[Proof of Theorem \ref{UPPER.BOUND.1}]
It follows from our assumption that for all $n > N(\varepsilon)$,
\[
\frac{p_{n+1}-p_n}{p_n p_{n+1}}\leq c(\varepsilon) \left( \frac{1}{n(n+1)} \right)^{1-\varepsilon}.
\]
Therefore, applying Rosser's Theorem, for sufficiently large $n$
\begin{align*}
p_{n+1} - p_n &\leq c(\varepsilon) \cdot p_n p_{n+1} \cdot \brackets{\frac{1}{n(n+1)}}^{1-\varepsilon} \\
& \leq 2c(\varepsilon) \cdot n^2 (\log n)^2 \cdot \brackets{\frac{1}{n}}^{2-2\varepsilon} \\
& \leq 2c(\varepsilon) \cdot (n^\varepsilon \log n)^2
\end{align*}
as required.
\end{proof}

Theorem \ref{UPPER.BOUND.2} provides a converse to the above, and requires a little more to prove.

\begin{proof}[Proof of Theorem \ref{UPPER.BOUND.2}]
Applying Rosser's Theorem, for sufficiently large $n$, 
\[
p_n \geq n \log n.
\]
 Fix $\varepsilon \geq 0$.  It follows from this estimate and our assumption that for sufficiently large $n$
\begin{align*}
\absolutevalues{\frac{1}{p_n} - \frac{1}{p_{n+1}}} & = \frac{p_{n+1} - p_n}{p_np_{n+1}} \\
& \leq  c(\varepsilon) \cdot \frac{ (n^\varepsilon \log n )^2}{n \cdot (n+1) \cdot \log n \cdot \log (n+1)} \\
& \leq c(\varepsilon) \cdot \brackets{\frac{1}{n(n+1)}}^{1-\varepsilon} \\
& = c(\varepsilon) \absolutevalues{\frac{1}{n} - \frac{1}{n+1}}^{1-\varepsilon}.
\end{align*}

Moreover, we can ``upgrade'' this estimate for $n < m < 2n$ using  the  concave  version of Jensen's inequality.  We obtain
\begin{align*}
\absolutevalues{\frac{1}{p_n} - \frac{1}{p_m}} & = \absolutevalues{\frac{1}{p_n} - \frac{1}{p_{n+1}}} + \dots  +\absolutevalues{\frac{1}{p_{m-1}} - \frac{1}{p_m}}  \\
& = \brackets{m-n} \cdot \frac{1}{m-n} \cdot \brackets{\absolutevalues{\frac{1}{p_n} - \frac{1}{p_{n+1}}} + \dots  +\absolutevalues{\frac{1}{p_{m-1}} - \frac{1}{p_m}}} \\
& \leq c(\varepsilon) \brackets{m-n} \cdot \frac{1}{m-n} \cdot \brackets{\brackets{\frac{1}{n} - \frac{1}{n+1}}^{1-\varepsilon} + \dots  +\brackets{\frac{1}{m-1} - \frac{1}{m}}^{1-\varepsilon}} \\
& \leq c(\varepsilon) \brackets{m-n} \cdot \brackets{\frac{1}{m-n}}^{1-\varepsilon} \brackets{\frac{1}{n} - \frac{1}{m}}^{1-\varepsilon} \\
& = c(\varepsilon) \brackets{m-n}^\varepsilon \cdot \brackets{\frac{1}{n} - \frac{1}{m}}^\varepsilon \cdot \brackets{\frac{1}{n} - \frac{1}{m}}^{1-2\varepsilon} \\
& \leq c(\varepsilon) \brackets{\frac{1}{n} - \frac{1}{m}}^{1-2\varepsilon}
\end{align*}
completing the proof.
\end{proof}

\noindent  {\large \bf Acknowledgements.}

H. Chen is thankful for  the excellent atmosphere for research provided by the University of St Andrews. The research of H. Chen was funded by China Scholarship Council (File No. 201906150102). J. M. Fraser was financially supported by an  EPSRC Standard Grant (EP/R015104/1) and a Leverhulme Trust Research Project Grant (RPG-2019-034).  The authors thank an anonymous referee for making helpful comments.

\end{document}